\documentclass[12pt]{article}
\usepackage{amssymb,latexsym,amsmath}

\newtheorem{thm}{Theorem}[section]

\newtheorem{cor}[thm]{Corollary}
\newtheorem{defn}{Definition}[section]

\def\Bx{\hfill $\Box$}

\begin{document}

\vskip 20pt
\begin{center}
{\bf REFINED RESTRICTED PERMUTATIONS AVOIDING SUBSETS
OF PATTERNS OF LENGTH THREE}
\vskip 15pt
{\bf Toufik Mansour}\\
{\it LaBRI, Universit\'e Bordeaux 1,
              351 cours de la Lib\'eration
              33405 Talence Cedex, France\\
        {\tt toufik@labri.fr} }

{\bf Aaron Robertson\footnote{
Homepage:  {\tt http://math.colgate.edu/$\sim$aaron/}
}
}\\
{\it Department of Mathematics,}
{\it Colgate University,
Hamilton, NY 13346}\\
{\tt aaron@math.colgate.edu}

\end{center}

\begin{abstract}
Define $S_n^k(T)$ to be the set of permutations of $\{1,2,\dots,n\}$
with exactly $k$ fixed points
which avoid all patterns in $T \subseteq S_m$.
We enumerate  $S_n^k(T)$, $T \subseteq S_3$, for all $|T| \geq 2$ 
and $0 \leq k \leq n$.
\end{abstract}

\vskip 30pt
\section{Introduction}
Let $\pi \in S_n$ be a permutation of $\{1,2,\dots,n\}$ written
in one-line notation.
Let $\alpha \in S_m$.
We say that $\pi$ {\it contains the pattern $\alpha$} if there exist
indices $i_1,i_2,\dots,i_m$ such that $\pi_{i_1} \pi_{i_2} \dots
\pi_{i_m}$ is equivalent to $\alpha$, where we define equivalence as
follows. Define $\overline{\pi}_{i_j}=
|\{x:\pi_{i_x} \leq \pi_{i_j}, 1 \leq x \leq m\}|$.  If
$\alpha = \overline{\pi}_{i_1} \overline{\pi}_{i_2} \dots
\overline{\pi}_{i_m}$ then
we say that $\alpha$ and $\pi_{i_1} \pi_{i_2} \dots \pi_{i_m}$ are
equivalent.  

Define $S_n^k(T)$, $T \subseteq S_m$, to be
the set of permutations with $k$ fixed
points which avoid all pattern in $T$.

In recent paper [RSZ], the authors
began the study of $S_n^k(T)$ for $T \subseteq S_3$, $|T|=1$.
In this paper we enumerate $S_n^k(T)$ for all
$T \subseteq S_3$ with $|T| \geq 2$.  Thus,  all that remains
 to enumerate are $S_n^k(123)$ and $S_n^k(231)$ for
all $0 \leq k \leq n$.

We start by giving two bijections which preserve
the parameter ``number of fixed points."
Let $I:S_n \rightarrow S_n$ be the group theoretical inverse.  Let
$R:S_n \rightarrow S_n$ be the reversal bijection,
i.e. $R(\pi_1 \pi_2 \dots \pi_n) = \pi_n \pi_{n-1} \dots \pi_1$, and
let $C:S_n \rightarrow S_n$ be the complement
bijection, i.e. $C(\pi_1 \pi_2 \dots \pi_n)=
(n+1-\pi_1) (n+1-\pi_2) \dots (n+1-\pi_n)$.

It is easy to see that the bijections
$I$ and $RC$  preserve
the number of fixed points in a permutation.
This reduces the number of cases that
we have to consider.

In the sequel, we will be using the
following notation.
Let $I_S$ be the characteristic function,
i.e. $I_S=1$ if $S$ is true and
$I_S=0$ if $S$ is false.
We also define $s_n^j(T)=0$ for $j<0$ or $j>n$
for any $T \subseteq S_3$.

\section{The Case $|T|=2$}

Using $I$ and $RC$ we see that we have the following
cases.

(1) $\overline{\{123,321\}}=\{\{123,321\}\}$

(2) $\overline{\{123,132\}}=\{\{123,132\},\{123,213\}\}$

(3) $\overline{\{123,231\}}=\{\{123,231\},\{123,312\}\}$

(4) $\overline{\{132,213\}}=\{\{132,213\}\}$

(5) $\overline{\{132,231\}}=\{\{132,231\},\{132,312\}\}$

(6) $\overline{\{132,321\}}=\{\{132,321\},\{213,321\}\}$

(7) $\overline{\{213,231\}}=\{\{213,231\},\{213,312\}\}$

(8) $\overline{\{231,312\}}=\{\{231,312\}\}$

(9) $\overline{\{231,321\}}=\{\{231,321\},\{312,321\}\}$

\begin{thm} $\{s_n^0(123,321)\}_{n \geq 0} = 1,0,1,2,4,0,0,\dots$,
$\{s_n^1(123,321)\}_{n \geq 0} =0,1,0,2,0,0,\dots$,
$\{s_n^2(123,321)\}_{n \geq 0}=0,0,1,0,0\dots$, and
$s_n^k(123,321)=0$ for all $3 \leq k \leq n$.
\end{thm}

{\bf Proof.} Obvious.
\Bx

\begin{thm} For $n \geq 1$, $s_n^k(123,132)=0$
for $3 \leq k \leq n$,
$$
\begin{array}{lc}
s_{2n+i}^2(123,132)=
\left\{
\begin{array}{ll}
\frac{4^{n-1}+2}{3}&\mathit{if \,\,}i=0\\\\
0&\mathit{if \,\,} i=1
\end{array}
\right., 
&
s_{2n+i}^1(123,132)=
\left\{
\begin{array}{ll}
\frac{2(4^{n-1}-1)}{3}&\mathit{if \,\,} i=0\\\\
\frac{4^{n}+2}{3}&\mathit{if \,\,}i=1
\end{array}
\right., and
\\\\
s_{2n+i}^0(123,132)=
\left\{
\begin{array}{ll}
4^{n-1}&\mathit{if \,\,} i=0\\\\
\frac{2(4^{n}-1)}{3}&\mathit{if \,\,}i=1
\end{array}
\right..
\end{array}
$$
\end{thm}

{\bf Proof.}  Let $\pi \in S_n(123,132)$. It is easy to prove by induction
that we must have
$\pi=(\pi(1), \pi(2), \dots,\pi(m))$ where
$\pi(j)=(t_{j-1}-1,t_{j-1}-2,\dots,t_j+1,t_{j-1})$
for some $n=t_0>t_1>t_2>\cdots>t_m=0$.
We call $\pi(j)$ a $j$-block of $\pi$.
(See [Man] for more details.)

Cleary $s_n^k(123,132)=0$ for
$3 \leq k \leq n$ since three fixed points
yield a $123$ pattern.  Hence,
we need only consider $k=0,1,2$.
However, using a result in [SS] and
the fact that $s_n^k(123,132)=0$ for
$3 \leq k \leq n$, we
immediately obtain the stated equation
for $s_n^0(123,231)$ upon proving
the formulas for $k=1,2$.  Hence, we
need only consider $k=1,2$.

We start with $k=2$.
Consider the {\it graph of $\pi$}, which consists
of the points  $(i,\pi_i)\in\mathbb{Z}^2$.
By the symmetries of the graph of $\pi$, we
see that the number of intersections with
the line $y=x$ (i.e. the number
of fixed points of $\pi$) is two if and only if
$m$ is odd
for $\pi=(\pi(1), \pi(2), \dots,\pi(m))$ 
and the intersections are in the $\frac{m-1}{2}$-block
of $\pi$.  So, the number of
elements in the $\frac{m-1}{2}$-block must
be even, as well as the total number of elements
of $\pi$.  Hence, we have
$s_n^2(123,132)=0$ for $n$ odd.
For $n$ even, let $2d$ be the number of elements
in the $\frac{m-1}{2}$-block.  Since the
other blocks cannot contain a fixed point,
we have $(s_{n-d}(123,132))^2$ permutations
for $d=1,2,\dots,n$.
Using a result in [SS], this gives
$1+\sum_{d=1}^{n-1} (2^{n-d-1})^2$ permutations.
Hence, $s_{2n}^2(123,132) =\frac{4^{n-1}+2}{3}$.

Next, we have $k=1$ and again we consider the
graph of $\pi$.  

For $n$ odd, we have a situation
similar to the $k=2$ case.  We must have the middle
block containing an odd number of elements.  Let 
this middle block have $2d+1$ elements.  
Using a result in [SS], this gives
$1+\sum_{d=1}^{n-1} (2^{n-d-1})^2$ permutations.
Hence, $s_{2n+1}^1(123,132) =\frac{4^{n-1}+2}{3}$.

For $n$ even, we have three possibilities.

1. $\pi(m)$ contains $2d_1$
elements, for some $2 \leq d_1 \leq n-1$, with one fixed point, or

2. $\pi(1)$ contain $2d_2$
of elements, for some $2 \leq d_2 \leq n-1$, with one fixed point, or

3. There exists exactly one block with $2d_3+1$ elements,
for some $2 \leq d_3 \leq n-1$, 
containing one fixed point. 

Letting $N=2n$, using a result in [SS], and summing over these three cases
we get
$$
\sum_{d_1=2}^{n-1} 2^{N-2d_1-1} + \sum_{d_2=2}^{n-1}
2^{N-2d_2-1}
+\sum_{d_3=1}^{n-1}(2^{n-d_3-1})^2
$$
permutations.
This gives $s_{2n}^1(123,132)=\frac{2(4^{n-1}-1)}{3}$,
thereby completing the proof.
\Bx

\begin{thm} For $n \geq 2$, $s_n^k(123,231)=0$
for $3 \leq k \leq n$,
$$
\begin{array}{l}
s_n^2(123,231)= 
\left\{
\begin{array}{ll}
\frac{n(n-6)}{24}+\frac{n}{2}&\mathit{if \,\,} n \equiv 0\, (
\mathit{mod\,\,}6)\\ \\
\frac{(n-1)(n+1)}{24}&\mathit{if \,\,} n \equiv 1,5 \, (
\mathit{mod\,\,}6)\\ \\
\frac{(n-4)(n-2)}{24}+\frac{n}{2}&\mathit{if \,\,} n \equiv 2,4 \, (
\mathit{mod\,\,}6)\\ \\
\frac{(n-3)(n+3)}{24}&\mathit{if \,\,} n \equiv 3 \, (
\mathit{mod\,\,}6)\\ \\
\end{array}
\right. ,
\\\\
s_n^1(123,231)=
\left\{
\begin{array}{ll}
\frac{n(n-6)}{12}+6{\frac{n+6}{6} \choose 2}&\mathit{if \,\,}
n
\equiv 0\, (
\mathit{mod\,\,}6)\\ \\
\frac{(n-3)(n-1)}{8}+\frac{(n-7)(n-1)}{12}
+ 6{\frac{n+5}{6} \choose 2}+\frac{n+2}{3}&\mathit{if \,\,} n
\equiv 1
\, (
\mathit{mod\,\,}6)\\ \\
\frac{n(n-2)}{12}+6{\frac{n+4}{6} \choose 2}&\mathit{if \,\,} n \equiv 2
\, (
\mathit{mod\,\,}6)\\ \\
\frac{(n-3)(n-1)}{8}+\frac{(n-5)(n-3)}{12}
+ 6{\frac{n+3}{6} \choose 2}+\frac{2n+3}{3}&\mathit{if \,\,} n \equiv 3 \,
(
\mathit{mod\,\,}6)\\ \\
\frac{(n-12)(n+2)}{12}+6{\frac{n+8}{6} \choose
2}&\mathit{if
\,\,} n
\equiv 4 \, (
\mathit{mod\,\,}6)\\ \\
\frac{(n-3)(n-1)}{8}+\frac{(n-5)(n-3)}{12}
+ 6{\frac{n+1}{6} \choose 2}+n&\mathit{if \,\,} n \equiv 5 \,
(
\mathit{mod\,\,}6)\\ \\
\end{array}
\right. , \,\, and
\end{array}
$$
and $s_n^0(123,231)={n \choose 2}+1 - s_n^1(123,231)
-s_n^2(123,231)$.
\end{thm}

{\bf Proof.}  Cleary $s_n^k(123,231)=0$ for
$3 \leq k \leq n$ since three fixed points
yield a $123$ pattern.  Hence,
we need only consider $k=0,1,2$.
However, using a result in [SS] and
the fact that $s_n^k(123,231)=0$ for
$3 \leq k \leq n$, we
immediately obtain the stated equation
for $s_n^0(123,231)$ upon proving
the formulas for $k=1,2$.  Hence, we
need only consider $k=1,2$.
   
Let  $\pi \in S_n(123,231)$ with
$\pi_i=n$.  If $i \neq 1$ then we must have
$\pi=(i-1)(i-2)\cdots 2 1 n (n-1) (n-2) \cdots i$
for any $2 \leq i \leq n$.  This case contributes
$n-1$ permutations to $S_n^1(123,231)$ for
$n$ odd and $\frac{n}{2}$ to 
$S_n^2(123,231)$ for $n$ even.  For $i=1$, we have
$\pi=n (n-1) \cdots (n-x+1) y (y-1) \cdots 1
(n-x) (n-x-1) \cdots (y+1)$ for some
$1 \leq x \leq n-2$ and $1 \leq y \leq n-x-1$,
or we have $\pi=n (n-1) \cdots 2 1$.

We now consider
the cases $k=1,2$ for $i=1$. 

We have three potential intervals in which fixed
points may reside:

(A) $n (n-1) \cdots (n-x+1)$, 

(B) $y (y-1) \cdots 1$,
or

(C) $(n-x) (n-x-1) \cdots (y+1)$

Hence, we have the
following situations in which a fixed point
appears.

A.  For $n$ odd we get one fixed point for any
$\frac{n+1}{2} \leq x \leq n-2$ and
$1 \leq y \leq n-x-1$, or if
$\pi=n (n-1) \cdots 1$.  The fixed point occurs
at position $\frac{n+1}{2}$.

B.  For $x+y$ odd we get one fixed point for any
$1 \leq x \leq \frac{n-2}{2}$ and
$x+1 \leq y \leq n-x-1$.  The fixed point
occurs at position $\frac{x+y+1}{2}$.

C.  For $n+y$ odd we get one fixed point for any
$1 \leq x \leq \frac{n-2}{2}$ and
$1 \leq y \leq n-2x-1$.  The fixed point
occurs at position $\frac{n+y+1}{2}$.

We now consider the case $k=2$ ($i=1$).  To have two distinct
fixed points we must have both 
$B$ and $C$ both occuring.  Hence,
we get
$$
\sum_{x=1}^{\lfloor \frac{n-2}{3} \rfloor} 
\sum_{{{y=x+1 \atop n+y \,\, odd} \atop x+y \,\, odd} }^{n-2x-1} 1 
= 
\left\{
\begin{array}{ll}
\frac{n(n-6)}{24}&\mathit{if} \,\, n \equiv 0 \, (
\mathit{mod\,\,}6)\\ \\
\frac{(n-1)(n+1)}{24}&\mathit{if} \,\, n \equiv 1,5 \, (
\mathit{mod\,\,}6)\\ \\
\frac{(n-4)(n-2)}{24}&\mathit{if} \,\, n \equiv 2,4 \, (
\mathit{mod\,\,}6)\\ \\
\frac{(n-3)(n+3)}{24}&\mathit{if} \,\, n \equiv 3 \, (
\mathit{mod\,\,}6)\\ \\
\end{array}
\right. 
$$
permutations with two fixed points in this case.

We now consider the case $k=1$ ($i=1$).  We may have only A
occur, only B occur, or only C occur:

Only A occurs:  In this case we require $n$ to be odd
and we get
$$
1+\sum_{x=\frac{n+1}{2}}^{n-2} \sum_{y=1}^{n-x-1} 1
= \frac{(n-3)(n-1)}{8} +1
$$
permutations with one fixed point.

Only B occurs:  In this case we get
$$
\sum_{x=1}^{\lfloor \frac{n-1}{3} \rfloor}
\sum_{{y=n-2x \atop x+y \,\, odd}}^{n-x-1} \!\! 1 + \!\!\!
\sum_{x=\lfloor \frac{n-1}{3} \rfloor+1}^{\lfloor \frac{n-2}{2}
\rfloor} \sum_{{y=x \atop x+y \,\, odd}}^{n-x-1} \!\! 1
+ \sum_{{x=1 \atop x \,\, odd}}^{\lfloor \frac{n-2}{3} \rfloor}
\sum_{{y=x+1 \atop y \,\, even}}^{n-2x-1} \!\! 1 
=\left\{
\begin{array}{ll}
\frac{n(n-6)}{24}+3{\frac{n+6}{6} \choose 2}&\mathit{if \,\,}
n
\equiv 0\, (
\mathit{mod\,\,}6)\\ \\
\frac{(n-15)(n-1)}{24}
+ 3{\frac{n+5}{6} \choose 2}&\mathit{if \,\,} n
\equiv 1
\, (
\mathit{mod\,\,}6)\\ \\
\frac{n(n-2)}{24}+3{\frac{n+4}{6} \choose 2}&\mathit{if \,\,} n \equiv 2
\, (
\mathit{mod\,\,}6)\\ \\
\frac{(n-9)(n-3)}{24}
+ 3{\frac{n+3}{6} \choose 2}&\mathit{if \,\,} n \equiv 3 \,
(
\mathit{mod\,\,}6)\\ \\
\frac{(n-12)(n+2)}{24}+3{\frac{n+8}{6} \choose
2}&\mathit{if
\,\,} n
\equiv 4 \, (
\mathit{mod\,\,}6)\\ \\
\frac{(n-5)(n-3)}{24}
+ 3{\frac{n+1}{6} \choose 2}&\mathit{if \,\,} n \equiv 5 \,
(
\mathit{mod\,\,}6)\\ \\
\end{array}
\right. 
$$
permutations with one fixed point.

Only C occurs:  In this case we get
$$
\sum_{x=1}^{\lfloor \frac{n-1}{3} \rfloor}
\sum_{{y=1 \atop n+y \,\, odd}}^x  \!\!1  + \!\!\!
\sum_{x=\lfloor \frac{n-1}{3} \rfloor+1}^{\lfloor \frac{n-2}{2}
\rfloor} \sum_{{{y=1 \atop n+y \,\, odd} \atop x+y \,\, even}}^{n-2x-1}
\!\! 1
+ \sum_{x=1}^{\lfloor \frac{n-2}{3} \rfloor}
\sum_{{y=x+1 \atop n+y \,\, odd}}^{n-2x-1} \!\! 1 
=\left\{
\begin{array}{ll}
\frac{n(n-6)}{24}+3{\frac{n+6}{6} \choose 2}&\mathit{if \,\,}
n
\equiv 0\, (
\mathit{mod\,\,}6)\\ \\
\frac{(n-15)(n-1)}{24}
+ 3{\frac{n+5}{6} \choose 2}&\mathit{if \,\,} n
\equiv 1
\, (
\mathit{mod\,\,}6)\\ \\
\frac{n(n-2)}{24}+3{\frac{n+4}{6} \choose 2}&\mathit{if \,\,} n \equiv 2
\, (
\mathit{mod\,\,}6)\\ \\
\frac{(n-9)(n-3)}{24}
+ 3{\frac{n+3}{6} \choose 2}&\mathit{if \,\,} n \equiv 3 \,
(
\mathit{mod\,\,}6)\\ \\
\frac{(n-12)(n+2)}{24}+3{\frac{n+8}{6} \choose
2}&\mathit{if
\,\,} n
\equiv 4 \, (
\mathit{mod\,\,}6)\\ \\
\frac{(n-5)(n-3)}{24}
+ 3{\frac{n+1}{6} \choose 2}&\mathit{if \,\,} n \equiv 5 \,
(
\mathit{mod\,\,}6)\\ \\
\end{array}
\right. 
$$
permutations with one fixed point.

Summing over all cases yields the stated formula.
\Bx

\begin{thm} For $n \geq 1$,
$s_n^{n-1}(213,132)=0$, $s_n^n(213,132)=1$,
$$
s_{2n+i}^0(213,132)=
\left\{
\begin{array}{ll}
\frac{5 \cdot 4^{n-1}-2}{3}&\mathit{if \,\,} i=0\\\\
\frac{2(4^{n}-1)}{3}&\mathit{if \,\,}i=1
\end{array}
\right.,
$$
and, for $1\leq k \leq n-1$,
$$
\begin{array}{l}
s_{2n+i}^{2k}(213,132)=
\left\{
\begin{array}{ll}
{4^{n-k-1}}&\mathit{if \,\,}i=0\\\\
0&\mathit{if \,\,} i=1
\end{array}
\right., and
\\\\
s_{2n+i}^{2k+1}(213,132)=
\left\{
\begin{array}{ll}
0&\mathit{if \,\,} i=0\\\\
{4^{n-k-1}}&\mathit{if \,\,}i=1
\end{array}
\right..
\end{array}
$$

\end{thm}

{\bf Proof.}  The proof is very similar to the
proof of Theorem 2.2 and we
provide a sketch of this proof.

It is easy to prove by induction that
$\pi=(\pi(1), \pi(2), \dots,\pi(m))$ where
$\pi(j)=(t_{j}+1,t_{j}+2,\dots,t_{j-1}-1,t_{j-1})$
for some $0=t_0<t_1<t_2<\cdots<t_m=n$.
We call $\pi(j)$ a $j$-block of $\pi$.
(See [Man] for more details.)

Recall that the graph of
$\pi$ consists of the points $(i,\pi_i) \in \mathbb{Z}^2$.

For $2k$ fixed points, we see that the middle block of
$\pi=(\pi(1),\dots,\pi(m))$ must be of length $2k$
and must contain all $2k$
fixed points and that $n$ must be even.  Hence, we have
$s_{2n}^{2k} = (s_{n-k}(213,132))^2=4^{n-k-1}$ 
and  $s_{2n+1}^{2k}=0$ for $1 \leq k \leq n-1$.

Similarly, for $2k+1$ fixed points, we have
$s_{2n+1}^{2k+1}=(s_{n-k}(213,132))^2=4^{n-k-1}$ 
and  $s_{2n}^{2k+1}=0$ for $1 \leq k \leq n-1$.

Using a result in [SS] gives the stated formula
for $s_n^0(213,132)$.
\Bx

\begin{thm} For $n \geq 3$, $s_n^0(132,231)=\frac13(2^{n-1}+(-1)^n)$,
$s_n^k(132,231)=\frac23(2^{n-k}+(-1)^{n-k+1})$
for $1 \leq k \leq n-2$, $s_n^{n-1}(132,231)=0$, and 
$s_n^n(132,231)=1$.  
\end{thm}

{\bf Proof.}  Let $\pi \in S_n(132,231)$.  Clearly,
we must have either $\pi_1=n$ or $\pi_n=n$.
For the case $\pi_n=n$ we have
$s_{n-1}^{k-1}(132,231)$ permutations.  Hence, we need
only consider $\pi_1=n$.  Write
$\pi=(n,\pi(1),1,\pi(2))$.  We must have
$\pi(1)$ be a decreasing sequence and $\pi(2)$ be
an increasing sequence.  Hence, we cannot have
a fixed point in $\pi(2)$ and we can have at most
one fixed point in $\pi(1)$.  Thus, for $k \geq 2$
the case $\pi_1=n$ is empty giving
$s_n^k(132,231)=s_{n-1}^{k-1}(132,231)$ for $k \geq 2$.

It remains to consider $k=0,1$ for
$\pi=(n,\pi(1),1,\pi(2))$.  We consider $k=0$
first (the case $k=1$ will follow immediately).
We may place the elements $2,3,\dots,n-1$ in $2^{n-2}$
ways such that $\pi(1)$ is decreasing and $\pi(2)$
is increasing.  To see this, insert the elements
$2,3,\dots,n-1$ in reverse order at $1$'s position 
and slide the element all the way
to the left or all the way to right.  However, we must subtract those
permutations which contain a fixed point, which
may only occur in $\pi(1)$.

Let $f$ be a fixed point and
write $\pi(1) = (\sigma(1),f,\sigma(2))$.  We have
${n-f-1 \choose f-2}$ choices for the elements
in $\sigma(1)$.  Once these are chosen, there
are then $2^{f-2}$ ways to place the
elements $2,3,\dots,f-1$ into $\sigma(2)$ and
$\pi(2)$ such that $\sigma(2)$ is decreasing and
$\pi(2)$ is increasing.  Thus, we have
$$
s_n^0(132,231) = 2^{n-2} - \sum_{f=2}^{\lfloor 
\frac{n}{2} \rfloor} {n-f-1 \choose f-2}2^{f-1}.
$$
From here, it is easy to check that
$s_n^0(132,231)=s_{n-1}^0(132,231)+2s_{n-2}^0(132,231)$,
which gives $s_n^0(132,231)=\frac13(2^{n-1}+(-1)^n)$
for $s_1^0(132,231)=0$ and $ s_2^0(132,231)=1$.

For the case $k=1$, the above argument shows that
$$
s_n^1(132,231) = \sum_{f=2}^{\lfloor 
\frac{n}{2} \rfloor} {n-f-1 \choose f-2}2^{f-1}.
$$
From here, it is easy to check that
$s_n^1(132,231)=s_{n-1}^1(132,231)+2s_{n-2}^1(132,231)$,
which gives $s_n^1(132,231)=\frac23(2^{n-1}+(-1)^n)$
for $s_1^1(132,231)=1$ and $ s_2^1(132,231)=0$.
\Bx

{\it Remark.}  $s_n^0(132,231)=J_{n-2}$ where
$J_n$ is the $n^\mathit{th}$ Jacobsthal number.

\begin{thm} For $n \geq 1$, $s_n^k(132,321)=n-k-1$ for 
$0 \leq k \leq n-1$ and
$s_n^n(132,321)=1$.
\end{thm}

{\bf Proof.}  Let $\pi \in S_n(132,321)$.
It is easy to see that in order to avoid
both $132$ and $321$ we have either
$\pi=(\pi',n)$ where $\pi' \in S_{n-1}(132,321)$
or $\pi=(j+1) (j+2) \cdots (n-1) n 1 2 \dots j$ for some
$ 1\leq j \leq n-1$.

Adjusting for the number of fixed points we see that
$$s_n^k(132,321)=s_{n-1}^{k-1}(132,321) + \sum_{j=1}^{n-1}I_{k=0}.$$

A straightforward induction on $k$ for $0 \leq k<n$ 
finishes the proof.
\Bx

As a corollary of Theorem 2.6 we rederive a result given
in [SS].

\begin{cor} $s_n(132,321)={n \choose 2} +1$
\end{cor}

{\bf Proof.}  Summing over $0 \leq k \leq n$ we obtain
the result immediately.

\begin{thm} For $n \geq 3$,
$s_n^0(213,231)=\frac13(2^{n-1}+(-1)^n)$,
$s_n^k(213,231)=\frac23(2^{n-k}+(-1)^{n-k+1})$
for $1 \leq k \leq n-2$, $s_n^{n-1}(213,231)=0$,
and $s_n^n(213,231)=1$. 
\end{thm}  

{\bf Proof.}
Using the fact that
$s_n^k(213,231)=s_n^k(213,312)$ and using the patterns
$213$ and $312$, the proof of this theorem is nearly 
identical to the proof of Theorem 2.5 and is
left to the reader.
\Bx

We now find the need for the following
definition.

\begin{defn} Let $S_1,S_2 \subseteq S_m$.  If $s_n^k(S_1)=s_n^k(S_2)$ for
all $0 \leq k \leq n$ for all $n\geq m$ we say that $S_1$ and
$S_2$  are Super-Wilf equivalent. 
\end{defn}

{\it Remark.}  Using the above definition, Theorems 1.5 and 1.8, and
results from [RSZ] we see that
$321$, $132$, and $213$ are Super-Wilf equivalent,
$231$ and $312$ are Super-Wilf equivalent, and
$\{132,231\}$, $\{132,312\}$, $\{213,231\}$ and $\{213,312\}$ are
Super-Wilf equivalent.

\begin{thm} For $n \geq 1$,
$$
s_n^k(231,312)=
\left\{
\begin{array}{ll}
2^{\frac{n-k-2}{2}}\left( {\frac{n+k}{2} \choose 
\frac{n-k}{2}} +
{{\frac{n+k-2}{2}} \choose 
\frac{n-k}{2}} \right)&\mathit{for \,\,}n+k \,\, \mathit{even}\\ \\
0&\mathit{otherwise}\\
\end{array}
\right.
$$

\end{thm}

{\bf Proof.}  For $\pi \in S_n(231,312)$ it is easy
to see that we must have $\pi=\pi(1) n (n-1) \cdots j$
with $\pi(1) \in S_{j-1}^k(231,312)$ for
some $1 \leq j \leq n$.
Hence, we have
$$
s_n^k(231,312) = 
\sum_{ {{j=1} \atop {n-j \,\, odd}} }^{n-1}
s_{j-1}^k(231,312) + \sum_{{ {j=1} \atop { {n-j \,\, even} \atop {or
\,\,n=j}}  }}^{n} s_{j-1}^{k-1}(231,312).
$$

From here we see that
$s_n^k(231,312)=2s_{n-2}^k(231,312) + s_{n-1}^{k-1}(231,312)$.
Hence, using initial conditions it is easy to check
that the formula given is correct.
\Bx

\begin{thm} Let $G_k(x)$ be the generating
function for $\{s_n^k(231,321)\}_{n \geq 0}$.
Then $G_k(x) = \frac{x^k(1-x)^{k+1}}{(1-x-x^2)^{k+1}}$.
In particular, $s_n^0(231,321)=F_{n-2}$,
for $n \geq 2$, where
$\{F_n\}_{n \geq 0}$ is the Fibonacci sequence
initialized by $F_0=F_1=1$.

\end{thm}

{\bf Proof.}  Let $\pi \in S_n(231,321)$.  It is easy to
see that we must have $\pi=(\pi',\sigma)$
where $\pi \in S_{n-j}^{k - I_{j=1}}(231,321)$
and $\sigma=n (n-j+1) (n-j+2) \cdots (n-1)$ for
some $1 \leq j \leq n$.

Hence, we have
$s_n^k(231,321)=s_{n-1}^{k-1}(231,321) + \sum_{j=2}^n s_{n-j}^k(231,321)$
for $0 \leq k \leq n$.  This implies that
$s_n^k(231,321)=s_{n-1}^{k-1}(231,321)+s_{n-2}^k(231,321)
+s_{n-1}^k(231,321)-s_{n-2}^{k-1}(231,321)$.

For $k=0$ we have
$s_n^0(231,321)=s_{n-1}^0(231,321)+s_{n-2}^0(231,321)
$, from which it follows that
$G_0(x) = \frac{1-x}{1-x-x^2}$.
For $k>0$ we have
$G_k(x)=\frac{x(1-x)}{1-x-x^2}G_{k-1}(x)$.
From here all of the results follow.
\Bx

The following corollary rederives a result in [SS] .

\begin{cor}
$\sum_{n \geq 0} s_n(231,321) x^n = \frac{1-x}{1-2x}$,
i.e. $s_n(231,321)=2^{n-1}$ for $n \geq 1$.
\end{cor}

{\bf Proof.}  Summing $\sum_{k \geq 0} G_k(x)$ gives
the desired result.
\Bx

\section{The Case $|T|=3$}

Using $I$ and $RC$ we see that we have the following
cases, except that case (1) is grouped together
because both $123$ and $321$ are to be avoided,
which is not possible for $n \geq 5$.

(1) $\overline{\{123,132,321\}}=\{\{123,132,321\},
\{123,213,321\},\{123,231,321\},\{123,312,321\}\}$

(2) $\overline{\{123,132,213\}}=\{\{123,132,213\}\}$

(3) $\overline{\{123,132,231\}}=\{\{123,132,231\},\{123,132,312\},
\{123,213,231\},\{123,213,312\}\}$

(4) $\overline{\{123,231,312\}}=\{\{123,231,312\}\}$

(5) $\overline{\{132,213,231\}}=\{\{132,213,231\},\{132,213,312\}\}$

(6) $\overline{\{132,213,321\}}=\{\{132,213,321\}\}$

(7) $\overline{\{132,231,312\}}=\{\{132,231,312\},\{213,231,312\}\}$

(8)
$\overline{\{132,231,321\}}=\{\{132,231,321\},\{132,312,321\},
\{213,231,321\},\{213,312,321\}\}$

(9) $\overline{\{231,312,321\}}=\{\{231,312,321\}\}$

\begin{thm} Let $\alpha \in \{132,213,231,312\}$. 
Then
$$
\begin{array}{l}
\{s_n^0(123,321)\}_{n
\geq 0} = 1,0,1,2,1,0,0,\dots
\,\, \mathit{for \,\,} \alpha \neq 231,312\\
\{s_n^0(123,321)\}_{n
\geq 0} = 1,0,1,1,1,0,0,\dots
\,\, \mathit{for \,\,} \alpha =231\,\, \mathit{or \,\,} \alpha=312\\
\{s_n^1(123,\alpha,321)\}_{n \geq 0} =0,1,0,2,0,0,\dots
\,\, \mathit{for \,\,} \alpha \neq 132,213\\
\{s_n^1(123,\alpha,321)\}_{n \geq 0} =0,1,0,1,0,0,\dots
\,\, \mathit{for \,\,} \alpha =132\,\, \mathit{or \,\,} \alpha=213\\
\{s_n^2(123,\alpha,321)\}_{n \geq 0}=0,0,1,0,0\dots,\,\,
\mathit{for \,\, any \,\,} \alpha, \,\,
\mathit{and}\\
\end{array}
$$
$s_n^k(123,\alpha,321)=0$ for any $\alpha$, for all $3 \leq
k\leq n$.
\end{thm}

{\bf Proof.} Obvious.
\Bx

\begin{thm} For $n \geq 3$, $s_n^k(123,132,213)=0$ for
$3 \leq k \leq n$, $s_n^2(123,132,213)=F^2_{\frac{n-2}{2}} I_{n \,\,
even}$, $s_n^1(123,132,213)=F^2_{\frac{n-1}{2}} I_{n \,\, odd}$,
and
$$
s_n^0(123,132,213) =
\left\{
\begin{array}{ll}
F_n-F^2_{\frac{n-2}{2}}&\mathit{if \,\,} n \,\, \mathit{is \,\, even}\\\\
F_n-F^2_{\frac{n-1}{2}}&\mathit{if \,\,} n \,\, \mathit{is \,\, even}\\
\end{array}
\right.,
$$
where $\{F_n\}_{n \geq 0}$ is the Fibonacci sequence
intialized by $F_0=F_1=1$.
\end{thm}

{\bf Proof.}  Clearly, $s_n^k(123,132,213)=0$ for $k \geq 3$ since
three fixed points are a $123$ pattern.   We may then use a
result from
from [SS] to get
$s_n^0(123,132,213)=s_n(123,132,213)-s_n^1(123,132,213)
-s_n^2(123,132,213)=F_n-s_n^1(123,132,213)
-s_n^2(123,132,213)$. 
Hence, we only
need consider $k=1,2$.  

We start with $k=2$.  The fixed points must be adjacent in order
to avoid all of $123, 132, 213$.  Let the fixed points be
$j$ and $j+1$.  Write $\pi=(\pi(1),j,j+1,\pi(2))$.  Then
$\pi(1)$ is on the elements $\{n-j+2,\dots,n\}$.  To avoid
the $123$ pattern, we require $n-j+2 \geq j+2$ and
$n-j+1 \leq j+1$.  Hence, $j=\frac{n}{2}$ (so that $n$ must
be even).  Thus, since it is not possible to have fixed
points in $\pi(1)$ or $\pi(2)$ we have
$s_n^2(123,132,213)=(s_{\frac{n-2}{2}}(123,132,213))^2 I_{n \,\, even}$,
which, using a result from [SS], gives the stated formula.

Now consider $k=1$.  Let $j$ be the fixed point.  Note that
$j \neq 1,n$ in order to avoid all of $123, 132, 213$.
Write $\pi=(\pi(1),j,\pi(2))$.  In order to avoid
$123$ and $213$, there exists at most one $x \in \pi(1)$
with $x<j$.  Assume, for a contradiction, that such an $x$ exists. Then we
must have
$x=\pi_{j-1}$ in order to avoid the $132$ pattern.  Next, for
$z \in \pi(2)$ we must have $z<j$ else $xjz$ is a $123$
pattern.  This forces $x=j-1$ to be a fixed point, a contradiction.
Hence, no such $x$ exists so that  for all $i \in \pi(1)$
we have $i>j$.  In order to avoid
$123$ and $132$, there exists at most one $y \in \pi(2)$
with $y>j$.  Assume, for a contradiction, that such a $y$ exists.  Then we
must have $y=\pi_{j+1}$ in order to avoid the $213$ pattern.  Furthermore,
for all $i \in \pi(1)$ we must have $i>j$ so as not to
have a $213$ pattern with $ijy$.  This forces
$y=j+1$ to be a fixed point, a contradiction.  Hence, no
such $y$ exists.  Since we have $i \in \pi(1)$, $i>j$ and
$k \in \pi(2)$, $k<j$, and for $j$ to be a fixed point
we require $n$ to be odd so that $j=\frac{n+1}{2}$, we get
$s_n^1(123,132,213)=(s_{\frac{n-1}{2}}(123,132,213))^2I_{n \,\, odd}$,
which, using a result from [SS], gives the stated formula.
\Bx

\begin{thm} For $n \geq 3$, $s_n^0(123,132,231)=\lfloor \frac{n}{2}
\rfloor$, 
$s_n^1(123,132,231)=\lfloor \frac{n}{2} \rfloor +(-1)^{n+1}$,
$s_n^2(123,132,231)=\frac{1}{2} (1+(-1)^n)$,
and $s_n^k(123,132,231)=0$ for $3 \leq k \leq n$.
\end{thm}

{\bf Proof.}  Let $\pi \in S_n(123,132,231)$.
It is easy to see that we must have
$\pi=n (n-1) \cdots (n-j+1) (n-j-1) (n-j-2) \cdots 2 1 (n-j)$
for some $0 \leq j \leq n-1$.  From here the results
follow easily.
\Bx

\begin{thm} For $n \geq 3$, $s_n^0(123,231,312)=s_n^2(123,231,312)
=\frac{n}{2}(1-I_{n \,\, odd})$,\\ $s_n^1(123,231,312)
=n(1-I_{n \,\, even})$, and
$s_n^k(123,231,312)=0$ for
$3 \leq k \leq n$.
\end{thm}

{\bf Proof.}  Let $\pi \in S_n(123,231,312)$.  It is
easy to see that we must have
$\pi = j (j-1) \cdots 1 n (n-1) \cdots (j+1)$.  From here
the results follow easily.
\Bx

\begin{thm} For $n \geq 3$, $s_n^0(132,213,231)=\lfloor \frac{n}{2}
\rfloor + (\frac{n}{2}+1)I_{n \,\, even}$, $s_n^1(132,213,231)
=  \lfloor \frac{n}{2} \rfloor I_{n \,\, odd}$,
and $s_n^k(132,213,231)=I_{n=k}$ for 
$2 \leq k \leq n$.
\end{thm}

{\bf Proof.}  Let $\pi \in S_n(132,213,231)$.  It is
easy to see that we must have
$\pi = n(n-1) \cdots (n-j+1) 1 2 \cdots (n-j)$
for some $1 \leq j \leq n$.  From here
the results follow easily.
\Bx

\begin{thm} For $n \geq 3$, $s_n^0(132,213,321)=n-1$ and
$s_n^k(132,213,321)=I_{n=k}$ for $1 \leq k \leq n$.
\end{thm}

{\bf Proof.}  Let $\pi \in S_n(132,213,321)$.  Then
we must have
$\pi = j (j+1) \cdots n 1 2 \cdots (j-1)$ for
some $1 \leq j \leq n$.  From here
the results follow easily.
\Bx

\begin{thm} For $n \geq 2$, $s_n^0(132,231,312)=\frac{1}{2}(1+(-1)^n)$,
$$s_n^{2k}(132,231,312)=
\left\{
\begin{array}{lll}
1+(-1)^n&\hskip 10pt&\mathit{if \,\,} n>2k\\
1&\hskip 10pt&\mathit{if \,\,} n=2k\\
\end{array}
\right., and
$$
$$s_n^{2k-1}(132,231,312)=
\left\{
\begin{array}{lll}
1+(-1)^{n+1}&\hskip 10pt&\mathit{if \,\,} n>2k-1\\
1&\hskip 10pt&\mathit{if \,\,} n=2k-1\\
\end{array}
\right..
$$
\end{thm}

{\bf Proof.}  Let $\pi \in S_n(132,231,312)$.  It is
easy to see that we must have
$\pi = j (j-1) \cdots 1 (j+1)(j+2) \cdots n$ for
some $1 \leq j \leq n$.  From here
the results follow easily.
\Bx

\begin{thm}  For $n \geq 3$,
$$s_n^{k}(132,231,321)=
\left\{
\begin{array}{lll}
1&\hskip 10pt&\mathit{if \,\,} 0 \leq k \leq n-2\\
0&\hskip 10pt&\mathit{if \,\,} k=n-1\\
1&\hskip 10pt&\mathit{if \,\,} k=n\\
\end{array}
\right..
$$
\end{thm}

{\bf Proof.}  Let $\pi \in S_n(132,231,321)$.  It is
easy to see that we must have
$\pi = j 1 2 \cdots (j-1)(j+1) \cdots n$ for
some $1 \leq j \leq n$.  From here
the results follow easily.
\Bx

\begin{thm} For  $n \geq 3$ and
$0 \leq k \leq n$, $s_n^k(231,312,321)=\binom{ \frac{n+k}{2}
}{k}I_{n +k \,\, even}$. 
\end{thm}

{\bf Proof.}  Let $\pi \in S_n^k(231,312,321)$.
Then $\pi$ must be of the form
$(1,\pi')$ or $(2,1,\pi'')$ where $\pi' \in S_{n-1}^{k-1}(231,312,321)$
and $\pi'' \in S_{n-1}^k(231,312,321)$.
This gives us $s_n^k(231,312,321) = s_{n-1}^{k-1} (231,312,321)
+ s_{n-2}^k(231,312,321)$.  A straightforward induction
on $n+k$ finishes the proof.
\Bx

\section{The Cases $|T| \geq 4$}

The cases $|T| \geq 4$ are easy and in fact
$s_n^k(T) \in \{0,1,2\}$ for all $T \subseteq S_3$, $|T| \geq 4$,
for all $n \geq 1$ and $0 \leq k \leq n$.

\section*{References}
[Man] T. Mansour, Permutations Avoiding a
Pattern from $S_k$ and at Least Two Patterns
from $S_3$, {\it Ars Combinatoria} {\bf 62}
(2001),

[RSZ] A. Robertson, D. Saracino, and D. Zeilberger,
Refined Restricted Permutations, 
arXiv: math.CO/0203033

[SS] R. Simion and F. Schmidt, Restricted Permutations,
{\it European Journal of Combinatorics} {\bf 6} (1985), 383-406.

\end{document}